%Paper: alg-geom/9506008
%From: Steven Bradlow <bradlow@vortex.math.uiuc.edu>
%Date: Mon, 5 Jun 95 17:50:03 -0600

%%%%%%%%%%%%%%%%%%%%%%%%%%%%%%%%%%%%%%%
%%  Hermitian-Einstein Inequalities  %%
%% and Harder-Narasimhan Filtrations %%
%%%%%%%%%%%%%%%%%%%%%%%%%%%%%%%%%%%%%%%
%%        Steven B. Bradlow          %%
%%%%%%%%%%%%%%%%%%%%%%%%%%%%%%%%%%%%%%%
%%             AMSTeX                %%
%%%%%%%%%%%%%%%%%%%%%%%%%%%%%%%%%%%%%%%
%%%%%%%%%%%%%%%%%%%%%%%%%%%%%%%%%%
%  THIS FILE LAST MODIFIED ON:   %
%         March 14, 1995         %
%%%%%%%%%%%%%%%%%%%%%%%%%%%%%%%%%%

\input amstex
\documentstyle{amsppt}
%\refstyle {A}
%\widestnumber\key{BGP}
\magnification=1200
\NoBlackBoxes
%%%%%%%%%%%%%%%%%%%%%%%%%%%%%%%%%%%
%general
%%%%%%%%%%%%%%%%%%%%%%%%%%%%%%%%%%%
\def \-->{\longrightarrow}
\def \kler{K\"ahler}

\def \tstable{$\tau$-stable}
\def \tstability{$\tau$-stability}

\def \curvdh{F_{\dbare,H}}

\def \dbar{\overline{\partial}}
\def \dbare{\overline{\partial}_E}

\def \veqtn{\frac{\sqrt{1}}{2\pi}\Lambda\curvdh+
        \frac{1}{2\pi}\phi\otimes\phi^*=\tau\bold I}

\def\today{\ifcase\month\or
  January\or February\or March\or April\or May\or June\or
  July\or August\or September\or October\or November\or
  December\fi
  \space\number\day, \number\year}

\topmatter
\title Hermitian-Einstein Inequalities and Harder-Narasimhan
Filtrations\endtitle
\author Steven B. Bradlow \endauthor
\address Department of Mathematics, University of Illinois,
Urbana, IL 61801 \endaddress
\email  bradlow\@uiuc.edu \endemail
\thanks Supported in part by NSF grant DMS-9303545 \endthanks
\abstract
Unstable holomorphic bundles can be described algebraically
by Harder-Narasimhan filtrations. In this note we
show how such filtrations correspond
to the existence of special metrics defined
by Hermitian-Einstein inequalities.
\endabstract
\endtopmatter
%\hfil\today
\bigskip

\document

%%%%%%%%%%%%%%%%%%%%%%%%%%%%%%%%%%%%%%%%%%%%%%%%%%%%%%%%%
\heading 1. Introduction \endheading
%%%%%%%%%%%%%%%%%%%%%%%%%%%%%%%%%%%%%%%%%%%%%%%%%%%%%%%%%

By a theorem of Uhlenbeck and Yau (cf. [10]), the stability of a holomorphic
bundle
over a closed \kler\ manifold can be detected by the existence of a Hermitian
bundle
metric which satisfies the Hermitian-Einstein equation.
In [6], Guan showed how a modification of this equation can be used to
quantify deviations from stability for non-stable
bundles.  Such failure of stability can also be measured
algebraically, namely by means of a Harder-Narasimhan filtration. These two
methods
naturally invite comparison, and that is the primary goal of this
note.  We also give some conditions on
non-stable bundles which are sufficient for the existence of solutions to the
modified Hermitian-Einstein equations.

Let $(X,\omega)$\ be a closed \kler\ manifold, and let $\Cal E\-->X$\ be a
rank
R holomorphic bundle over $X$. Denote the underlying smooth bundle by $E$.
Then
$\Cal E$\ corresponds to $E$\ together with an integrable $\dbar$-operator,
which we denote by $\dbare$. To avoid unnecessary extra notation, we will use
$\Cal E$\ to refer both to the bundle and the corresponding locally free
coherent analytic sheaf.  We define the
slope of any coherent analytic subsheaf $\Cal E'\subset \Cal E$\ in
the usual way, i.e. by
$$\mu(\Cal E') = \frac{\int_X {c_1(\Cal E')\wedge\omega^{n-1}}}{rank(\Cal E')}
.$$
The bundle is stable (respectively semistable) if for all
coherent subsheaves with $0<rank(\Cal E')<R$\ we have
$$\mu(\Cal E')<\mu(\Cal E)\  (\mu(\Cal E')\le\mu(\Cal E)\ ) .$$

The Hermitian-Einstein equation is an equation for an
Hermitian bundle metric. It is given in terms of the
curvature of the unique connection compatible with both the metric
and the bundle's holomorphic structure.  If we denote the metric
by $H$, and the curvature of the associated connection by
$\curvdh$, then the equation reads

$$\frac{\sqrt{-1}}{2\pi}\Lambda \curvdh = \mu(\Cal E)\bold I .$$

Here $\bold I$\ is the identity bundle automorphism, and
$\Lambda \curvdh$\ denotes the bundle endomorphism obtained by
taking the contraction of $\curvdh$\ against the \kler\ form.
Thus $\Lambda \curvdh\wedge\omega^n = \curvdh\wedge\omega^{n-1}$.

\proclaim{Theorem [10]}An indecomposable holomorphic bundle $\Cal E$\ is
stable
if and only if it supports a Hermitian metric satisfying the
Hermitian-Einstein equation.
\endproclaim

In the modification considered by Guan, a parameter is
introduced in the right hand side of the Hermitian-Einstein
equation.  Because of the Chern-Weil homomorphism , which
relates $\frac{\sqrt{-1}}{2\pi}\Lambda \curvdh$\ to the first
Chern class of $E$, there can be no solutions to the \it
equality \rm $\frac{\sqrt{-1}}{2\pi}\Lambda \curvdh = m\bold I$\
unless $m=\mu(\Cal E)$.  This topological constraint will however
permit solutions if the equality is replaced by an
inequality.  Notice that the expressions on both sides of such a condition are
Hermitian bundle endomorphisms. One can make sense of an inequality between
two
such endomorphisms by interpreting $A\le B$\ to mean that $A-B$\ is negative
semi-definite (with similar meaning for the inequality $A\ge B$).  Guan's
result is

\proclaim{Theorem 1 [6]} Let $\Cal E\--> X$\ be a holomorphic
bundle over $X$.  Let $m$(resp. $m'$) be a real number and
suppose that $\Cal E$\ supports a metric $H$\ such that

$$\frac{\sqrt{-1}}{2\pi}\Lambda \curvdh \le m\bold I \ \
\biggr{(}\text{resp.}\
\frac{\sqrt{-1}}{2\pi}\Lambda \curvdh \ge m'\bold I\ \biggl{)} ,\tag
1$$
i.e. such that $\frac{\sqrt{-1}}{2\pi}\Lambda \curvdh - m\bold I\
\bigr{(}\text{resp.}\  m'\bold I-\frac{\sqrt{-1}}{2\pi}\Lambda \curvdh\
\bigl{)}$\ is a negative semidefinite (Hermitian) bundle endomorphism.  Then
$$\mu(\Cal E')\le m\ \  \biggr{(}\text{resp.}\  \mu(\Cal E/\Cal E')
\ge m' \ \biggl{)}\tag 2$$
for all subsheaves $\Cal E'\subset \Cal E$.
\endproclaim

\noindent\bf{Remark:}\rm\  If $\mu(\Cal E)\ne 0$, then one can write
$m=t\mu(\Cal E),\ m'=t'\mu(\Cal E)$, which is the way the results are
presented
in [6].

A natural question to consider is for which
values (if any) of $m$\ or $m'$\ the equations in (1) have a solution.  That
is, if
$$\Cal M = \{m\ :\ \frac{\sqrt{-1}}{2\pi}\Lambda \curvdh\le
m\bold I\ \text{for some metric}\ H\ \}\ , $$

$$\Cal M' = \{m'\ :\ \frac{\sqrt{-1}}{2\pi}\Lambda \curvdh\ge
m'\bold I\ \text{for some metric}\ H\ \}\ , $$
what can one say about $\Cal M$\ and $\Cal M'$?

It is immediately clear that if non-empty, the sets $\Cal M $\  and $\Cal M'$\
are half-infinite intervals, with $\mu(\Cal E)$\ being a lower bound for $\Cal
M$ and an upper bound for $\Cal M'$ . In fact one can be a bit more precise,
even without any further work. We must first recall that for a holomorphic
bundle over a closed \kler\ manifold, the following is true (cf. [9])

\proclaim{Theorem 2 (Harder-Narasimhan Filtrations)}Given a holomorphic
bundle $\Cal E$\ over a closed \kler\ manifold $(X,\omega)$,
there is a unique filtration (called the Harder-Narasimhan
filtration) by subsheaves

$$0=\Cal E_0\subset\Cal E_1\subset\dots\Cal E_k=\Cal E\ ,\tag 3a$$
(where $\Cal E$\ is the sheaf associated to $E$)such that
$\Cal E_i/\Cal E_{i-1}$\ is the unique maximal semistable
subsheaf of $\Cal E/\Cal E_{i}$, for $1\le i\le k$.In
particular, the slope of the quotients are ordered such that

$$\mu(\Cal E_1)>\mu(\Cal E_2/\Cal E_1)>\dots\mu(\Cal
E_k/\Cal E_{k-1})\ .\tag 3b$$

If $\Cal E$\ is semistable, then there is a filtration by
subsheaves (called the Seshadri filtration)

$$0=\Cal E_0\subset\Cal E_1\subset\dots\Cal E_k=\Cal E\ ,\tag 3c$$
such that the quotients $\Cal E_i/\Cal E_{i-1}$\ are all
stable bundles and have slope $\mu(\Cal E_i/\Cal E_{i-1})=\mu(E)$.  Moreover,
$$Gr(\Cal E)=\Cal E_1\oplus\Cal E_2/\Cal E_1\oplus\dots\oplus\Cal
E/\Cal E_{k-1}\tag 3d$$
is uniquely determined up to an isomorphism.

If $X$ is a Riemann surface, then the terms in both of these
filtrations are locally free, i.e. are subbundles of $\Cal E$.
\endproclaim

\proclaim{Definition 3}Define  $\mu_1(\Cal E)$\ by
$$\mu_1(\Cal E)=\mu(\Cal E_1)\ ,$$
where $\Cal E_1\subset\Cal E$\ is the first term in the
Harder-Narasimhan filtration for $E$.
\endproclaim

Thus $\mu_1(\Cal E)$\ is the least upper bound for the slopes of all
subsheaves
$\Cal E'\subset \Cal E$.  It follows that $inf(\Cal M)\ge \mu_1(\Cal E)$.
Similarly, by using the correspondence between quotients of $\Cal E$\ and
subsheaves of the dual bundle $\Cal E^*$, one sees that $sup(\Cal M')\le
-\mu_1(\Cal E^*)$.

%%%%%%%%%%%%%%%%%%%%%%%%%%%%%%%%%%%%%%%%%%%%%%%%%%%%%%%%%
\heading 2. Over a compact Riemann surface \endheading
%%%%%%%%%%%%%%%%%%%%%%%%%%%%%%%%%%%%%%%%%%%%%%%%%%%%%%%%%

In the case where the base manifold $X$\ is a closed Riemann
surface, we can considerably strengthen the connection between the
Hermitian-Einstein inequalities and the Harder-Narasimhan filtrations. To be
precise:

\proclaim{Theorem 4} Let $\Cal E\--> X$\ be a holomorphic
bundle over a closed Riemann surface, and let $\Cal M$, $\Cal M'$\ be
as above.  Then the sets $\Cal M$\ and $\Cal M'$\ are non-empty, and
\roster
\item $inf(\Cal M)  = \mu_1(\Cal E)$,
\item $sup(\Cal M')  = -\mu_1(\Cal E^*)$.
\endroster
\endproclaim

The proof of Theorem 4 follows from a
more general result, namely:

\proclaim{Theorem 5}  Let $X$\ be a closed Riemann surface,
and suppose that the holomorphic bundle $\Cal E\--> X$\ has
Harder-Narasimhan filtration
$$0=\Cal E_0\subset\Cal E_1\subset \Cal E_2\subset \dots\subset \Cal E_k=\Cal
E\ .\tag 4$$
Let the ranks and slopes of $Q_i=\Cal E_{i}/\Cal E_{i-1}$\ be given by
$(r_i, \mu_i(\Cal E)),\ i=1,2,\dots,k$.

Then given any $\epsilon >0$\  there is a complex bundle automorphism, $g$ of
$E$, and an Hermitian metric
$H$\ on $E$\ such that
$$-\epsilon\bold I\le \frac{\sqrt{-1}}{2\pi}\Lambda F_{g(\dbare),H} -
 \pmatrix
\mu_1 \bold I_{r_1}& 0 & \hdots & 0\\
0 & \mu_2\bold I_{r_2} & \hdots & 0\\
\dots & \dots & \dots & \dots\\
0 & 0 & \hdots & \mu_k\bold I_{r_k}
\endpmatrix
\le \epsilon\bold I  $$

A solution corresponding to $\epsilon =0$\ can be found if
the smooth decomposition $\Cal E=\bigoplus_{i=1}^k
\Cal E_{i}/\Cal E_{i-1}$\ is a holomorphic decomposition, and each
quotient $Q_i=\Cal E_{i}/\Cal E_{i-1}$\ is a polystable
bundle.
\endproclaim

\demo{Proof}  Using Seshadri filtrations, we can refine the filtration for
$\Cal E$ to get a filtration in which all quotients are stable bundles. We may
thus assume that all the $Q_i$\ are stable bundles, and hence admit
Hermitian-Einstein metrics. Denote these by $K_i$, for $i=1,2,\dots,k$.  We
thus get a metric, say $K$, on the graded object $Gr(\Cal E)=\oplus Q_i$\ for
which
$$\frac{\sqrt{-1}}{2\pi}\Lambda F_{ \dbare^0,K}=
 \pmatrix
\mu_1 \bold I_{r_1}& 0 & \hdots & 0\\
0 & \mu_2\bold I_{r_2} & \hdots & 0\\
\dots & \dots & \dots & \dots\\
0 & 0 & \hdots & \mu_k\bold I_{r_k}
\endpmatrix\ , $$
where $\dbare^0$\ denotes the holomorphic structure on $Gr(\Cal E)$. We now
show how to construct a bundle automorphism $g$\ such that the curvature
$\frac{\sqrt{-1}}{2\pi}\Lambda F_{ g(\dbare),K}$\ is arbitrarily close to
$\frac{\sqrt{-1}}{2\pi}\Lambda F_{ \dbare^0,K}$.  Such arguments are well
known, and can be found in [1] (cf. section 8), and also [5].
The basic idea can be illustrated in the simplest non-trivial case, i.e. the
case where $\Cal E$\ is given as an extension of stable bundles, say
$$0\-->\Cal E_1\-->\Cal E\-->\Cal E_2\-->0\ .\tag 5$$
In this case we can pick Hermitian-Einstein metrics, $K_i$, on the bundles
$\Cal E_i$. Thus we have
$$\frac{\sqrt{-1}}{2\pi}\Lambda F_{\dbare^0,K} =
 \pmatrix
\mu_1 \bold I_{r_1}& 0 \\
0 &  \mu_2\bold I_{r_2}
\endpmatrix\ ,$$
where $\dbare^0$\ denotes the holomorphic structure on $\Cal E_1\oplus\Cal
E_2$, and $K=K_1\oplus K_2$.  Now with respect to the orthogonal splitting
determined by $K$, the holomorphic structure on $\Cal E$\ is given by
$$ \dbare=\pmatrix
\dbar_1& \beta \\
0 &  \dbar_2
\endpmatrix\ ,$$
where $\dbar_i$\ gives the holomorphic structure on $\Cal E_1$\ and
$\beta\in\Omega^{0,1}(Hom(E_2,E_1))$\ is the second fundamental form of the
inclusion $\Cal E_1\hookrightarrow\Cal E$. A straightforward computation gives
$$\frac{\sqrt{-1}}{2\pi}\Lambda F_{\dbare,K} =
\pmatrix
\mu_1 \bold I_{r_1} - \frac{\sqrt{-1}}{2\pi}\Lambda\beta\wedge\beta^*&
\frac{\sqrt{-1}}{2\pi}\Lambda d\beta\\
-\frac{\sqrt{-1}}{2\pi}\Lambda d\beta^* &
  \mu_2\bold I_{r_2}
-\frac{\sqrt{-1}}{2\pi}\Lambda \beta^*\wedge\beta
\endpmatrix\ ,$$
where $d$\ denotes covariant differentiation determined by the metric
connections on $\Cal E_1$\ and $\Cal E_2$.

Recall that the class $[\beta]\in H^1(Hom(E_2,E_1))$\ determines the
isomorphism class of (5) as an extension , while it is the corresponding point
in $\Bbb P(H^1(Hom(E_2,E_1)))$\ which gives the isomorphism class of the
bundle
$\Cal E$. Thus if we define
$$ \dbar_t=\pmatrix
\dbar_1& t\beta \\
0 &  \dbar_2
\endpmatrix\ ,$$
we get a  1-parameter family of extensions,  all of which are isomorphic to
$\Cal E$\ as bundles. In fact, the holomorphic structures $\dbar_t$\ and
$\dbar_1=\dbare$\ are related by the complex gauge transformation
$$g_t=\pmatrix
\bold I_{r_1}& 0 \\
0 &  t\bold I_{r_2}
\endpmatrix\ .$$
If we pick $\beta$\ to be the harmonic representative in its cohomology class,
then we find
$$\frac{\sqrt{-1}}{2\pi}\Lambda F_{\dbar_t,K}-
\pmatrix
\mu_1 \bold I_{r_1}& 0 \\
0 &  \mu_2\bold I_{r_2}
\endpmatrix\ = t^2
\pmatrix
 - \frac{\sqrt{-1}}{2\pi}\Lambda\beta\wedge\beta^*& 0\\
0 & -\frac{\sqrt{-1}}{2\pi}\Lambda \beta^*\wedge\beta
\endpmatrix\ .$$
The result now follows by taking small enough $t$.

Now, instead of assuming stability of $\Cal E_1$, suppose that the theorem
applies to $\Cal E_1$. That is, assume that $\Cal E_1$\ admits a metric such
that
$$-\epsilon\bold I\le \frac{\sqrt{-1}}{2\pi}\Lambda F_{\dbar_1,K_1} -
 \pmatrix
\mu^{(1)}_1 \bold I_{r_1}& 0 & \hdots & 0\\
0 & \mu^{(1)}_2\bold I_{r_2} & \hdots & 0\\
\dots & \dots & \dots & \dots\\
0 & 0 & \hdots & \mu^{(1)}_k\bold I_{r_k}
\endpmatrix
\le \epsilon\bold I \ ,\tag 6$$
where the $r_i$'s and $\mu^{(1)}_i$'s are the ranks and degrees of the
quotients in the filtration for $\Cal E_1$.  By exactly the same argument as
above, we see that we can find a 1-parameter family of complex bundle
automorphism, $g_t$, such that
$$\frac{\sqrt{-1}}{2\pi}\Lambda F_{\dbar_t,K}-
\pmatrix
\frac{\sqrt{-1}}{2\pi}\Lambda F_{\dbar_1,K_1}& 0 \\
0 &  \mu_2\bold I_{r_2}
\endpmatrix\
 = t^2
\pmatrix
 - \frac{\sqrt{-1}}{2\pi}\Lambda\beta\wedge\beta^*& 0\\
0 & -\frac{\sqrt{-1}}{2\pi}\Lambda \beta^*\wedge\beta
\endpmatrix\ .$$
Combining this with (6) shows that the theorem then applies to $\Cal E$. We
may
thus apply this method, one step at a time, to the filtration for $\Cal E$.
Since the gauge transformation at stage $j$\ leaves unchanged the filtration
up
to $\Cal E_j$, the composition of all the gauge transformations produces the
required result to prove the theorem.\qed
\enddemo

\bf Remark.\rm\ Another way to obtain this result is by considering the
gradient flow for the Yang-Mills functional on $\Cal C$,  the space of
holomorphic structure for $\Cal E$. This goes back to the work of [1], with
some important analytic details being supplied by [4]. Using a fixed
background
metric, $K$, we can define the functional $f:\Cal C\-->\Bbb R$\ by

$$f(\dbare)=||\frac{\sqrt{-1}}{2\pi}\Lambda F_{\dbare,K}-\mu(E)\bold
I||_{L^2}^2\ .$$

Via the identification of $\Cal C$\ and the space of unitary connections, this
is essentially the same
as the Yang-Mills functional $\Cal {YM}(D)=||F_D||^2_{L^2}$. The critical
points correspond to reducible holomorphic structures of the form $\Cal
E=\bigoplus \Cal E_i$, where each summand is stable. The results of [1]
(especially section 8) and [4] show that under the gradient flow, each
holomorphic bundle $\Cal E$ converges to a critical point corresponding
precisely to the graded object $Gr(\Cal E)$. One then uses the fact that the
gradient flow preserves the isomorphism class of $\Cal E$, and that at the
critical points $\frac{\sqrt{-1}}{2\pi}\Lambda F_{\dbare,K}$ has the required
diagonal form.

\demo{Proof of Theorem 4}
\it Part (1)\ \rm:  Recall that  $\mu_1(\Cal E)=\mu_1>\mu_2>\dots>\mu_k$.
Thus
$$\pmatrix
\mu_1\bold I_{n_1}& 0 & \hdots & 0\\
0 &\mu_2\bold I_{n_2} & \hdots & 0\\
\dots & \dots & \dots & \dots\\
0 & 0 & \hdots & \mu_k\bold I_{n_k}
\endpmatrix \le \mu_1(\Cal E)\bold I\ .$$
It then follows from Theorem 5 that for any$\epsilon >0$\  there is a complex
bundle automorphism, $g$ of $E$, and an Hermitian metric  $H$\ on $E$\ such
that
$$ \frac{\sqrt{-1}}{2\pi}\Lambda F_{g(\dbare),H}\le(\mu_1(\Cal
E)+\epsilon)\bold I\ .
\tag 4$$

Instead of using $g$\ to transform the $\dbar$-operator, we can use it to
change the metric on $E$.  If the new metric is $H_g=Hg^*g$, we have the
relation $$\frac{\sqrt{-1}}{2\pi}\Lambda
F_{g(\dbare),H}=g\circ\frac{\sqrt{-1}}{2\pi}\Lambda F_{\dbare,H_g}\circ
g^{-1}\
.$$
It follows then from (4) that
$$ \frac{\sqrt{-1}}{2\pi}\Lambda F_{\dbare,H_g}\le(\mu_1(\Cal
E)+\epsilon)\bold
I\ .$$
Together with the observation that $inf(\Cal M)\ge \mu_1(\Cal E)$, this proves
part (1)

\it Part (2)\ \rm: Given a connection $D$\ on $\Cal E$, there is a
corresponding induced connection $D^*$\ on $\Cal E^*$. Furthermore, if $H^*$\
is the Hermitian metric on $\Cal E^*$\ dual to the metric $H$\ on $\Cal E$,
then the metric connection  $D_{H^*}$\ on $\Cal E^*$\ is precisely the
connection corresponding to the metric connection $D_H$\ on $\Cal E$.
Suppose that with respect to a local holomorphic frame for $E$, the curvature
of $D_H$\ is $F_H$, and that with respect to the dual holomorphic frame for
$E^*$\ the curvature of $D_{H^*}$\ is $F_{H^*}$. Then
$$i\Lambda F_{H^*}=-i\Lambda F_H^t\ ,$$
where $F_H^t$\ denotes the transpose of $F_H$.  It follows that
$$\align
i\Lambda F_H \ge m'\bold I &\Longleftrightarrow -i\Lambda F_{H^*} \ge  m'\bold
I\\
&\Longleftrightarrow i\Lambda F _{H^*} \le -m'\bold I
\endalign$$
The result now follows from part (1). \hfill\qed
\enddemo

%%%%%%%%%%%%%%%%%%%%%%%%%%%%%%%%%%%%%%%%%%%%%%%%%%%%%%%%%%
\heading
3. Over higher dimensional \kler\ manifolds
\endheading
%%%%%%%%%%%%%%%%%%%%%%%%%%%%%%%%%%%%%%%%%%%%%%%%%%%%%%%%%%

For bundles over closed \kler\ manifolds
of dimension greater than two, the filtrations of holomorphic bundles are by
subsheaves, rather than subbundles.  The above arguments thus do not apply.
For related reasons, the proof of Theorem 4 based on the Morse theory of the
Yang-Mills functional also fails, with the method breaking down because of the
failure of the Yang-Mills gradient flow to converge. Such failure to converge
is caused by the ``bubbling" phenomenon on the space of connections.

In the case where $X$\ is not a Riemann surface, some information about
$inf(\Cal M)$\ and $sup(\Cal M')$\ can be obtained by relating the
Hermitian-Einstein
\it inequality \rm to the \it equality \rm given by the
$\tau$-Vortex equation. This is an equation for a metric on
a holomorphic bundle with a prescribed holomorphic section,
and has the form
$$\veqtn\ . $$
Here $\tau$\ is a real parameter which must lie in the range
$(\mu(\Cal E_1), \frac{R}{R-1}\mu(\Cal E))$\ (cf[2]).  Since
the bundle endomorphism $-\phi\otimes\phi^*$\ is non-positive,
it is apparent that

\proclaim{Lemma 7}Let $\Cal E\--> X$\ be a holomorphic bundle over
a closed \kler\ manifold of dimension $n\ge 1$. Fix the
parameter $\tau$.
Then there is a solution to the equation
$$\frac{\sqrt{-1}}{2\pi}\Lambda F_{H} \le \tau \bold I\
$$
if for some choice of section $\phi\in H^0(X,\Cal E)$\ there is a
solution to the $\tau$-Vortex equation
$$\veqtn\ .$$
\endproclaim

For a given pair $(\Cal E,\phi)$, there is a necessary and
sufficient condition for the $\tau$-vortex equation to have
a solution.  This is the \tstability condition for a
holomorphic pair.  We recall the result from [2]

\proclaim {Definition 8} Given a real number $\tau$, we say
that the pair $(E,\phi)$  is $\tau$-stable (resp.
 $\tau$-semistable) if the following
two conditions hold:
\roster
\item $\mu(\Cal E') <\tau$ (resp. $\leq \tau$),
for every holomorphic
subbundle $\Cal E_{\phi}\subset E$;
\item $\mu(\Cal E/\Cal E_{\phi})>\tau$ (resp. $\geq \tau$),
for every proper
holomorphic subbundle $\Cal E_{\phi}\subset \Cal E$ such that $\phi$
is a
section of  $\Cal E_{\phi}$.
\endroster
\endproclaim

\proclaim{Theorem 9}[2]
Suppose that $(\Cal E,\phi)$
is $\tau$-stable for a given value of the parameter
$\tau$. Then the $\tau$-Vortex equation
$$\frac{\sqrt{-1}}{2\pi}\Lambda\curvdh +\frac{1}{
2\pi}\phi\otimes\phi^\ast=\tau\bold I$$
considered as an equation for the metric $H$,
has a unique smooth solution.

Conversely, suppose that for a given value of $\tau$ there
is a Hermitian metric $H$ on $\Cal E$ such that the $\tau$-vortex
equation is satisfied on $(\Cal E,\phi)$.  Then $\Cal E$ splits
holomorphically as
$\Cal E=\Cal E_\phi\oplus \Cal E_s$, where
\roster
\item $\Cal E_s$, if not empty, is a direct sum of stable
bundles, each of slope $\tau\cdot \frac{Vol(X)}{4\pi}$,
\item $\Cal E_\phi$ contains the section $\phi$ and
$(\Cal E_\phi,\phi)$ is $\tau$-stable.
\endroster
\endproclaim

Notice that the split case $\Cal E=\Cal E_\phi\oplus \Cal E_s$
cannot occur unless $\tau\cdot \frac{Vol(X)}{4\pi}$
corresponds to the
slope of a subbundle, i.e. unless
$\tau\cdot \frac{Vol(X)}{4\pi}$\ is a rational number with
denominator less than the rank of $\Cal E$.  Hence,
for generic values of $\tau$ the summand $\Cal E_s$ is empty, and
$\tau$-stability is the necessary and sufficient condition
for  the
existence of solutions to the $\tau$-vortex equation.

For our present purposes it is convenient to define the
following
parameter .

\proclaim{Definition 10}Given a holomorphic pair $(\Cal E,\phi)$,
let
$$inf (\Cal E,\phi)=Min\{\mu(\Cal E/\Cal E_{\phi}\ :\ \ \Cal E_{\phi}\subset
\Cal E,\
\text{and}\ \phi\in H^0(X,\Cal E_{\phi})\}\ .$$
\endproclaim

We then have the result

\proclaim{Proposition 11}Let $(\Cal E,\phi)$\ be a holomorphic pair.
Let
$\mu_1(\Cal E)$\ be the slope of the first subbundle in the Harder-
Narasimhan
filtration for $E$.  Then  the pair is \tstable\ for some
value of $\tau$\
if and only
$$\mu_1(\Cal E) < inf(\Cal E,\phi) \ .$$
In that case, $\tau$\ lies in the interval $(\mu_1(\Cal E)\ ,\
inf(\Cal E,\phi))$.
\endproclaim

This gives us the corollary

\proclaim{Corollary 12}Let $\Cal E\--> X$\ be a holomorphic bundle
over a closed \kler\ manifold of dimension $n\ge 1$.
\roster
\item  If there is a section $\phi\in H^0(X,\Cal E)$\ such that
$\mu_1(\Cal E) < inf(\Cal E,\phi)$, then for all $m>\mu_1(\Cal E)$, there is a
solution to
$$\sqrt{-1}\Lambda\curvdh\le m\bold I\  .$$
\item If there is a section $\phi^*\in H^0(X,\Cal E^*)$\ such that
$\mu_1(\Cal E\*) < inf(\Cal E^*,\phi^*)$, then for all $m>\mu_1(\Cal E^*)$,
there is a solution to
$$\sqrt{-1}\Lambda\curvdh\ge -m\bold I\  .$$
\endroster
\endproclaim

\demo{Proof}
(1)  Given a section $\phi\in H^0(X,\Cal E)$\ such that $\mu_1 <
inf(\Cal E,\phi)$, the pair $(\Cal E,\phi)$\ is $\tau$-stable for any
$\mu_1<\tau<inf(\Cal E,\phi)$.  The result thus follows from Lemma 9.

(2)  Replace $\Cal E$\ by $\Cal E^*$\ in the proof of part (1)  \qed
\enddemo

%%%%%%%%%%%%%%%%%%%%%%%%%%%%%%%%%%%%%%%%%%%%%%%%%%%
%%%%       EXTRA REMARKS BEGIN HERE          %%%%%%
%%%%%%%%%%%%%%%%%%%%%%%%%%%%%%%%%%%%%%%%%%%%%%%%%%%

This result can be rephrased in an interesting way by using the interpretation
in [7], [8] of the vortex equations as a dimensional reduction of the
Hermitian-Einstein equations. As shown in [8], a holomorphic pair $(\Cal
E,\phi)$\ over $X$\ can be identified with a holomorphic extension over
$X\times\Bbb P^1$\ of the form
$$0\-->p^*\Cal E\-->\Bbb E\-->q^*\Cal O(2)\-->0\ .\tag 5$$
Here $p,q$\ are the projections from $X\times\Bbb P^1$\ onto $X$\ and $\Bbb
P^1$\ respectively, and $\Cal O(2)$\ is the degree two line bundle over $\Bbb
P^1$. The extension class of $\Bbb E$\ is related to a section  $ \phi\in
H^0(X,\Cal E)$\ by the isomorphism
$$\align
H^1(X\times\Bbb P^1,p^*\Cal E\otimes q^*\Cal O(-2))
&\cong  (H^0(X,\Cal E)\otimes H^1(\Bbb P^1,\Cal O(-2)))
\\
&\cong H^0(X,\Cal E)\ .
\endalign$$

To define stability or the Hermitian-Einstein equations on $\Bbb E$, we need
to
fix a \kler\ metric on $X\times\Bbb P^1$.  We consider the 1-parameter family
of such metrics corresponding to the \kler\ forms
$$\Omega_{\sigma}=p^*\omega+\sigma q^*\omega_{\Bbb P^1}\ .$$
Here $\omega$\ and $\omega_{\Bbb P^1}$ are the \kler\ forms on $X$\ and $\Bbb
P^1$\ respectively, and $\sigma$\ is a positive real parameter.  Notice that
there is a natural $SU(2)$\ action on $X\times\Bbb P^1$\ (which is trivial on
$X$). The \kler\ forms $\Omega_{\sigma}$\ are invariant under this action, and
there is a natural lift of the action to $\Bbb E$.  Garcia-Prada has shown

\proclaim {Theorem 13 [7] } Let $\sigma$\ and $\tau$\ be related by
$$\sigma=\frac{2Vol(X)}{(rank(\Cal E)+1)\tau - deg(\Cal E)}\ .\tag6$$
Under the above correspondence between holomorphic pairs on $X$\ and
holomorphic extensions on $X\times\Bbb P^1$, the following are equivalent:
\roster
\item There is a metric on $\Cal E$\ satisfying the $\tau$-vortex equation,
\item There is an $SU(2)$-invariant metric on $\Bbb E$\ satisfying the
Hermitian-Einstein equation with respect to $\Omega_{\sigma}$.
\endroster
\endproclaim

As a result of Theorem 9, and the correspondence between stability and
Hermitian-Einstein metrics (cf. [10]), this leads to the following

\proclaim{Theorem 14 [7]}Let $\sigma$\ and $\tau$\ be related as in Theorem
13.
Under the above correspondence between holomorphic pairs on $X$\ and
holomorphic extensions on $X\times\Bbb P^1$, the following are equivalent:
\roster
\item The pair $(\Cal E,\phi)$\ is $\tau$-stable,
\item The extension $\Bbb E$\ is stable with respect to $\Omega_{\sigma}$.
\endroster
\endproclaim

Corollary 12 (1) can thus be rephrased as

\proclaim {Corollary 15} Let $\Cal E\--> X$\ be a holomorphic bundle over a
closed \kler\ manifold of dimension $n\ge 1$. Suppose that there is a choice
of
$\sigma$\ and an extension
$$0\-->p^*\Cal E\-->\Bbb E\-->q^*\Cal O(2)\-->0\  $$
such that $\Bbb E$\ is stable with respect to $\Omega_{\sigma}$. Then for all
$m>\mu_1(\Cal E)$, there is a metric on $\Cal E$\ satisfying
$\sqrt{-1}\Lambda\curvdh\le m\bold I$.
\endproclaim
\demo{Proof}If  $\Bbb E$\ is stable with respect to $\Omega_{\sigma}$, then
the
corresponding pair $(\Cal E,\phi)$\ is $\tau$-stable, where $\tau$\ and
$\sigma$\ are related as in (6). The result now follows as before.\qed
\enddemo

A similar rephrasing of Corollary 12(2) is also possible.

%%%%%%%%%%%%%%%%%%%%%%%%%%%%%%%%%%%%%%%%%%%%%%%%%%%
%%%%       EXTRA REMARKS END HERE            %%%%%%
%%%%%%%%%%%%%%%%%%%%%%%%%%%%%%%%%%%%%%%%%%%%%%%%%%%

In the special case of rank two bundles, we can be even more
explicit.  We can use the fact that any section $\phi\in H^0(X,\Cal E)$\
generates a rank one subsheaf with torsion free quotient, and that this is the
only proper subsheaf which contains the section and has torsion free quotient.
Denoting this subsheaf by $[\phi]$, we thus get that
$$inf (\Cal E,\phi)=\mu(\Cal E/[\phi])=2\mu(\Cal E)-\mu([\phi])\ .$$
The holomorphic pair $(\Cal E,\phi)$\ will then be $\tau$-stable for some
value
of $\tau$\ if and only if $\mu([\phi])\ne \mu_1(\Cal E)$, i.e.
$\mu([\phi])<\mu_1(\Cal E)$.  Combining this with Corollary 12, we get

\proclaim{Corollary 16}Let $\Cal E\--> X$\ be a rank 2 holomorphic bundle over
a closed \kler\ manifold of dimension $n\ge 1$.
\roster
\item  Suppose that there is a section $\phi\in H^0(X,\Cal E)$\ such that
$deg([\phi])<\mu_1(\Cal E)$. Then for all $m>\mu_1(\Cal E)$, there is a
solution to
$$\sqrt{-1}\Lambda\curvdh\le m\bold I\ .$$
\item Suppose that there is a section $\phi^*\in H^0(X,\Cal E^*)$\ such that
$deg([\phi^*])<\mu_1(\Cal E^*)$. Then for all $m>\mu_1(\Cal E^*)$, there is a
solution to
$$\sqrt{-1}\Lambda\curvdh\ge -m\bold I\ .$$
\endroster
\endproclaim

\subheading{Remark}  The vortex equations and the definition of stable pairs
can be generalized in a way which conveniently takes into account the duality
between the two cases covered by Theorem 1. The new equations are coupled
equations for metrics on two bundles $\Cal E_1$\ and $\Cal E_2$\ over $X$, and
have the form
$$\frac{\sqrt{1}}{2\pi}\Lambda
F_{\dbar_1,H_1}+\frac{1}{2\pi}\Phi\otimes\Phi^*=\tau\bold I\ ,$$
$$\frac{\sqrt{1}}{2\pi}\Lambda
F_{\dbar_2,H_2}-\frac{1}{2\pi}\Phi^*\otimes\Phi=\tau'\bold I\ .$$
The section $\Phi$\ is now a section of $H^0(X,Hom(\Cal E_2,\Cal E_1))$, and
the constants $\tau$\ and $\tau'$\ are related by the constraint
$$r_1\tau + r_2\tau' = deg(E_1)+deg(E_2)\ ,$$
where $r_1$\ is the rank of $\Cal E_1$\ etc.

These equations where introduced in [7].  It is clear that the solutions to
these coupled vortex equations provide solutions to the inequalities
$\sqrt{-1}\Lambda F_{\dbar_1,H_1}\le \tau\bold I$\ on $\Cal E_1$, and
$\sqrt{-1}\Lambda F_{\dbar_2,H_2}\ge \tau'\bold I$\ on $\Cal E_2$.  In [7] and
[3] it is shown that the existence of such solutions is related to a stability
criterion (also called $\tau$-stability) for the triple $(\Cal E_1,\Cal
E_2,\Phi)$. The definition of $\tau$-stability is a slope condition on
subtriples, i.e. on $(\Cal E'_1,\Cal E'_2,\Phi')$\ where for $i=1,2$\ $\Cal
E'_i$\ is a  rank $r'_i$\ subsheaf of $\Cal E_i$\ and $\Phi'\in H^0(X,Hom(\Cal
E'_2,\Cal E'_1))$\ is such that the obvious diagram commutes (cf. [7]).  For
all such subtriples let
$$\theta_{\tau}(\Cal E'_1,\Cal E'_2)=(\mu({\Cal E}'_1\oplus {\Cal
E}'_2)-\tau)-\frac{r_2'}{r_2}\frac{r_1+r_2}{r_1'+r_2'}
 (\mu({\Cal E}_1\oplus {\Cal E}_2)-\tau)\ ,$$
and define the triple to be $\tau$-stable if $\theta_{\tau}(\Cal E'_1,\Cal
E'_2)<0$\ for all proper subtriples.  Then (cf. [3]) $\tau$-stability implies
the existence of metrics on $\Cal E_1$\ and $\Cal E_2$\ satisfying the coupled
vortex equations, and the following observation can then be made:

\proclaim {Proposition 14}Let $\Cal E\--> X$\ be a holomorphic bundle
over a closed \kler\ manifold of dimension $n\ge 1$.
\roster
\item Suppose there is a bundle $\Cal F\--> X$,  a holomorphic section
$\Phi\in
H^0(X, \Cal E\otimes\Cal F^*)$, and a real number $\tau$\ such that $(\Cal
E,\Cal F,\Phi)$\ is a $\tau$-stable triple. Then there is a solution to the
inequality
$\sqrt{-1}\Lambda \curvdh\le \tau\bold I$\ on $\Cal E$.
\item Suppose there is a bundle $\Cal F\--> X$,  a holomorphic section
$\Phi\in
H^0(X, \Cal F\otimes\Cal E^*)$, and a real number $\tau$\ such that $(\Cal
F,\Cal E,\Phi)$\ is a $\tau$-stable triple. Then there is a solution to the
inequality
$\sqrt{-1}\Lambda \curvdh \ge \tau'\bold I$\ on $\Cal E$, where $\tau$\ and
$\tau'$\ are related as above.
\endroster
\endproclaim

\enddocument
\end